\documentclass[12pt]{article}

\usepackage{epsf,psfrag,amsmath,amsfonts, amssymb}

\begin{document}
\bibliographystyle{abbrv}
\pagenumbering{arabic}
\raggedbottom

\title{A Central Limit Theorem \\ for non-overlapping return times}
\author{Oliver Johnson\thanks{Statistical Laboratory, 
University of Cambridge, Wilberforce Road, Cambridge CB3 0WB, UK. 
E-mail: \tt{otj1000@cam.ac.uk}}}
\date{\today}
\maketitle
\hfuzz20pt

\newtheorem{theorem}{Theorem}[section]
\newtheorem{lemma}[theorem]{Lemma}
\newtheorem{proposition}[theorem]{Proposition}
\newtheorem{remark}[theorem]{Remark}
\newtheorem{corollary}[theorem]{Corollary}
\newtheorem{conjecture}[theorem]{Conjecture}
\newtheorem{definition}[theorem]{Definition}
\newtheorem{example}[theorem]{Example}
\newtheorem{condition}{Condition}
\newtheorem{main}[theorem]{Theorem}
\renewcommand{\thefigure}{\thesection.\arabic{figure}}
\setlength{\parskip}{\parsep}
\setlength{\parindent}{0pt}
\setcounter{tocdepth}{1}

\def \outlineby #1#2#3{\vbox{\hrule\hbox{\vrule\kern #1%
\vbox{\kern #2 #3\kern #2}\kern #1\vrule}\hrule}}%
\def \endbox {\outlineby{4pt}{4pt}{}}%
\newenvironment{proof}
{\noindent{\bf Proof\ }}{{\hfill \endbox
}\par\vskip2\parsep}
\newenvironment{pfof}[2]{\removelastskip\vspace{6pt}\noindent
 {\it Proof  #1.}~\rm#2}{\par\vspace{6pt}}

\newcommand{\cond}{\stackrel{d}{\longrightarrow}}
\newcommand{\geo}[1]{{\rm Geom\,}(#1)}
\newcommand{\expo}[1]{{\rm Expo\,}(#1)}
\newcommand{\iv}{{\cal V}}
\newcommand{\A}{{\cal A}}
\newcommand{\vc}[1]{\mathbf{#1}}
\newcommand{\I}{\mathbb I}
\newcommand{\cov}{{\rm{Cov}}}
\newcommand{\ep}{{\mathbb E}}
\newcommand{\pr}{{\mathbb P}}
\newcommand{\zee}{{\mathbb Z}}
\newcommand{\var}{{\rm{Var}}}
\newcommand{\ta}{^{(\tau)}}
\newcommand{\tends}{\rightarrow \infty}
\newcommand{\blah}[1]{}

\newlength{\dsl}
\newcommand{\doublesub}[2]{\settowidth{\dsl}{$\scriptstyle
#1$}\parbox{\dsl}{\scriptsize\centering { \normalsize $\scriptstyle #1$}
\\ {\normalsize $\scriptstyle #2$}}}

\begin{abstract}
\noindent Define the non-overlapping 
return time of a random process to be the number 
of blocks that we wait before a particular block reappears. We prove a
Central Limit Theorem based on these  return times. 
This result has applications
to entropy estimation, and to the problem of determining if digits
have come from an independent equidistributed sequence. In the case of
an equidistributed sequence, we use an argument based on negative
association to prove convergence under weaker conditions.
\end{abstract}




\section{Introduction and main theorem}
\subsection{Statement of Problem}
Given a sample $(Z_1, \ldots Z_n)$ from a
random process taking values in an alphabet $\A$,
we would like to estimate the entropy of the process. In
general, this is a hard problem, though if the process is assumed to be
independent or stationary, some progress can be made.

In particular, given a sequence of binary bits, determining
whether the bits were generated by an independent equidistributed process
has applications to problems in
cryptography and number theory, as described in Section \ref{sec:applications}.

Our approach is as follows:
we first
partition the sample of $(Z_i)$ into blocks of size $\ell$. That is,
writing $Z_a^b$ for $(Z_a, Z_{a+1}, \ldots ,Z_b)$, we define block
random variables
$X_i = Z_{(i-1)\ell+1}^{i\ell}$, so
each $X_i \in \A^{\ell}$. Then, given the first $k$ blocks $X_1, \ldots,
X_k$, we count how long it takes for these blocks to reappear.
\begin{definition}[Non-overlapping return time] \label{def:defines}
For a given $k$, define random variable 
$$ S_j = \min \{ t \geq 1: X_{j+t} = X_j \}, \mbox{ for $j=1, \ldots, k$,}$$
to be the return time of the $j$th block.
\end{definition}
It appears that this definition dates back to Maurer \cite{maurer}.
The main result of this paper is that if the number and size of blocks
grow appropriately, then the $S_j$ satisfy a Central Limit Theorem:
\begin{theorem} \label{thm:main} 
Suppose that $(Z_i)$ is an independent identically distributed
finite alphabet process with entropy $H$. 
Write $q_{\max} < 1$ for the maximum
probability of any symbol.
If, as block length $\ell \tends$, the number of blocks $k(\ell) \tends$ in 
such a way that $\lim_{\ell \tends} k(\ell)^{3/2} \ell q_{\max}^{\ell} = 0$, 
then
$$ \frac{
\sum_{i=1}^{k(\ell)} (\log S_i - \ell H \log 2 + \gamma)}{\sqrt{
k(\ell) \pi^2/6  }} \cond
N(0,1).$$
\end{theorem}
\begin{remark} \mbox{  } 

\begin{enumerate} 
\item Here $\cond$ denotes convergence in distribution, and $\gamma$ is Euler's
constant $\gamma = 0.577216\ldots$. 
\item Note that to fit in with conventions in information theory, the 
entropy $H$ is calculated using logarithms to base 2. If entropy were 
calculated using natural logarithms, the $\log 2$ term could be omitted.
\item For equidistributed processes, where each symbol occurs with
probability $q$, Proposition \ref{prop:sclt} shows that
we can relax the assumption on the
rate of convergence of $k(\ell)$, to only require that
$k(\ell) \ell q^{\ell} \rightarrow 0$.
Indeed, simulations in Section \ref{sec:comp} suggest
that a weaker condition, namely $k(\ell) \ell 2^{-H \ell} \rightarrow 0$,
 may be sufficient to ensure convergence for all finite alphabet independent
processes with entropy $H$. The faster rate of convergence 
in the case of equidistribution is useful, since we will often test
a null hypothesis of equidistribution.
\end{enumerate}
\end{remark}

\begin{corollary}
Under the conditions of Theorem \ref{thm:main} above, the estimator
$$ \widehat{H} = 
\frac{ \sum_{i=1}^{k(\ell)} (\log S_i + \gamma)}{\ell \log 2
\sqrt{ k(\ell) \pi^2/6}}
\sim AN \left( H, \frac{1}{l^2} \right),$$
so is asymptotically normal and consistently estimates the entropy.
\end{corollary} 
If  the  process  is  independent  and  equidistributed  on  a  finite
alphabet,  each  block occurs  at  each  time  with probability  $p  =
|\A|^{-\ell}$. Hence each of the $S_j$ are geometric random variables,
with  $\pr(S_j  =  r)  =  p(1-p)^{r-1}$  (we  call  this  a  $\geo{p}$
variable).  In general, if  the process is independent and identically
distributed  (IID)  with  entropy  $H$, it  satisfies  the  Asymptotic
Equipartition Property (AEP), so that asymptotically there are $\simeq
2^{H  \ell}$   blocks  of  length  $\ell$  which   appear,  each  with
probability $\simeq  2^{-H \ell}$, so
conditioned  on  the value  of  $X_j$,  $S_j$  is a  geometric  random
variable.  However, even though the symbols $Z_i$ are independent, the
return  times  $S_j$ are  dependent,  so  we  need to  understand  the
dependence structure to prove Theorem \ref{thm:main}.

In Section \ref{sec:similar} we describe some results concerning similar 
return time definitions made by other authors. Section 
\ref{sec:applications} describes possible two applications of these results.
In Section \ref{sec:pfmain} we prove Theorem \ref{thm:main}, using an 
argument based on asymptotic independence.  We transform into a similar
problem, by eliminating the possibility of early matches.
Section \ref{sec:negass} gives
a proof under weaker conditions for equidistributed random
variables, by using negative association. Section \ref{sec:comp} contains
the results of some simulations.

In future work, we hope to extend these results to general stationary 
processes,
under a suitable mixing condition, and to prove similar results for other 
definitions of match length.
\subsection{Other similar definitions} \label{sec:similar}
We briefly describe some other work concerning similar quantities. 
This is by no means an exhaustive list, but merely gives a flavour of 
some of the alternative approaches which exist. 
\begin{enumerate}
\item{\label{item:overlap}
\begin{definition}[Overlapping return time] \label{def:definet}
Define random variables $T_k$ to be the time before the block $Z_1^k$ is next 
seen:
$$ T_k = \min \{ t \geq 1: Z_1^k = Z_{t+1}^{t+k} \}.$$
\end{definition}
Kac's Lemma \cite{kac} shows that for any stationary ergodic process
 $\ep[T_k | Z_1^k = z_1^k] = 1/\pr(Z_1^k = z_1^k)$. This was 
developed by Kim \cite{kim}, who gave the limiting behaviour of $\ep[T_k 
\pr(Z_1^k)]$ for independent processes, 
and by Wyner \cite{wyner4}, who showed that
\begin{theorem} \label{thm:wyner} If $(Z_i)$ is a Markov chain then
$$ \lim_{n \tends} \frac{ \log_2 T_n - n H}{\sqrt{n \iv}} \sim N(0,1),$$
\end{theorem}
Here $\iv$ is the information variance, 
$\lim_{n \tends} \var(-\log_2 \pr(Z_1^n))/n$.
For independent processes: \\ $\iv = \var(-\log_2 \pr(Z_1)) = 
\sum_i \pr(Z_1 = z_i) (-\log_2 \pr(Z_1 = z_i) - H)^2$.

See also Corollary 2 of Kontoyiannis \cite{kontoyiannis2}, who showed 
that this result holds
for general stationary $(Z_i)$ under explicit mixing conditions. Wyner and Ziv 
\cite{wyner}, Ornstein and Weiss \cite{ornstein} and Gao \cite{gao} 
consider similar quantities.}
\item{\label{item:grass}
\begin{definition}[Grassberger prefix]\label{def:match}
Given $n$ and $1\leq i \leq n$ define
$$ R_{i,n}(Z_1^n)= \inf\left\{ t: Z_i^{i+t-1}\neq Z_j^{j+t-1} \\ 
\mbox{ for all } j\neq i\right\}.$$
\end{definition}
In words, $R_{i,n}(Z_{1}^n)$
is the length of the shortest string started
at position $i$ different from all the others
of equal length started at $j$ for $1\leq j\leq n$.
This quantity was introduced by Grassberger \cite{grassberger}, and
studied by Kontoyiannis and Suhov
\cite{kontoyiannis}, Quas \cite{quas} and Shields \cite{shields},
\cite{shields2}, partly it allows
good entropy estimation for an ergodic process
with a suitable degree of mixing.
For example, Theorem 1 of \cite{kontoyiannis} shows:
\begin{theorem} \label{thm:kont}
If finite alphabet process $(Z_i)$ is ergodic with entropy $H$,
and satisfies a Doeblin condition then
$$ \lim_{n\to\infty} \frac{1}{n}
\sum_{i=1}^{n} \frac{\log n}{R_{i,n}(Z_1^n )} =H, \mbox{ almost surely}.$$
\end{theorem}
}
\item{ \label{item:lz}
{\bf Lempel-Ziv coding}
Another problem with similar features is 
finding the asymptotic behaviour of the number of codewords in
the Lempel--Ziv parsing (see Section 12.10 of Cover and Thomas \cite{cover}).
Ziv \cite{ziv3} made a conjecture concerning the number of
codewords. However, Aldous and Shields \cite{aldous} were only
able to resolve the problem for IID equidistributed processes, and it took
careful analysis by Jacquet and Szpankowski \cite{jacquet} to extend their
results to IID asymmetric processes. For example Theorem 1A of \cite{jacquet}
shows that
\begin{theorem} \label{thm:jacquet} 
Given an binary asymmetric ( $\pr(Z_1 = 0) \neq 1/2$)
IID process, then $L_m$,
the total length of the $m$ words in a Lempel--Ziv tree, satisfies
$$ \lim_{m \tends} \frac{ L_m - \ep L_m}{\sqrt{ \var(L_m)}} \sim N(0,1),$$
where $\ep L_m = (m \log_2 m)/H + O(m)$, $\var(L_m) = \iv m \log_2 m/H^3 
+ O(m)$.
\end{theorem}}
\end{enumerate}
Notice that these Theorems \ref{thm:wyner}, \ref{thm:kont} and 
\ref{thm:jacquet}
differ in character from our Theorem \ref{thm:main}. For example, Theorem
\ref{thm:kont} proves a law of large numbers for the Grassberger prefixes,
showing that a statistic based on them acts as an entropy estimator. However,
it does not tell us the rate of convergence of the estimator.
Similarly, although Theorems \ref{thm:wyner}
and \ref{thm:jacquet} give asymptotic normality, they both refer to statistics
calculated with respect to
one fixed point. It is possible that this fixed point
could be unrepresentative, and so our result is stronger in the sense that
it averages over a number of different starting points.
\subsection{Applications} \label{sec:applications}
We briefly describe two applications of these results.
\begin{enumerate}
\item{{\bf (Cryptography)}
The problem of deciding whether binary bits $(Z_i)$ were generated by an IID
equidistributed process arises in cryptography.
Bits generated in this way can be used as a perfectly
secure one--time pad to transmit a message $(Y_i)$. This system is secure 
in the sense that the 
transmitted bits $Y_i \oplus Z_i$ are independent of the message,
so no inference about $Y_i$ can be made from them. Equivalently
Shannon's Second Coding Theorem (see for example
Theorem 8.7.1 of \cite{cover}) implies that the binary 
symmetric channel with
error probability $p =1/2$ has capacity $C=0$. If the $(Z_i)$ were not
independent and equidistributed, given large enough $n$, it may be possible 
to infer
properties of the $(Z_i)$ and perhaps read the message $(Y_i)$.  }
\item{{\bf (Number Theory)} Recall that  a number is said to be normal
to base $b$  if the limiting proportion of each digit  in its base $b$
expansion  is $1/b$.  A number  which is  normal to  all bases  $b$ is
referred to as  normal.  Ergodic theory shows that  almost all numbers
are normal,  but it is  hard to prove  that any particular  number has
this property.  For example,  Bailey and Crandall  \cite{bailey} prove
that a  particular class of numbers (including  the so-called Stoneham
and Korobov numbers)  has the normal property.  On  the other hand, in
the  same paper  \cite{bailey} they  discuss the  fact  that constants
including  $\pi$, $e$,  $\ln 2$  and $\zeta(3)$  are not  known  to be
normal.   Weisstein \cite{weissteinnormal} gives  a review  of results
concerning normal  numbers.  An informal statement of  the property of
normality to  base $b$ is  that the digits  of the number `look  as if
they were generated by an  IID equidistributed process', which we hope
to be able to test.}
\end{enumerate}
Kim \cite{kim} gave computational results concerning the speed of convergence
of estimators based on overlapping matches, hoping to detect
processes which are not Bernoulli. Similarly
Bradley and Suhov \cite{bradley2} used theoretical results
concerning the Grassberger prefixes (see Definition \ref{def:match})
to consider the normality of constants such as $\pi$, $e$ and $\gamma$. 
We give some computational results in Section \ref{sec:comp}.
\section{Proof of Main Theorem} 
\label{sec:pfmain}
\subsection{Avoiding early matches}

The difficulty in analysing the dependence structure of the random variables 
$S_i$ introduced in 
Definition \ref{def:defines} is that `early matches' can occur at 
$i$. That is, it may be that
$S_i \leq k - i$. The possibility of early matches leads to a complicated 
situation of case splitting, according to where such early matches occur.

To avoid this, we introduce a very similar sequence of random 
variables $(R_j)$ in Definition
\ref{def:definer} below. We use two main ideas to prove Theorem \ref{thm:main},
 the Central Limit Theorem for the $S_i$. 

First $S_i = R_i + (k - i)$, unless there has been an early 
match. By controlling the probability of an early match, we 
show that a suitably scaled version of
$\sum \log S_i - \sum \log R_i$ tends to zero in probability, and so
a limit law for the $R_i$ passes over to a limit law for the $S_i$. 
The formal statement is
given in Lemma \ref{lem:roomcorr} below.

Then, in Proposition \ref{prop:charcl}, we establish
a Central Limit Theorem for the $R_i$, using explicit bounds on 
conditional probabilities, which show that the variables are asymptotically
independent.
\begin{definition} \label{def:definer} 
Given a realisation of $X_1, \ldots, X_k$, we define $D$,
the set of positions which do not see an early match. That is
$$ D = \{ i: X_i \neq X_{j} \mbox{ for $j = i+1, \ldots, k$} \}.$$
For each $i \in D$, we take $b_i = X_i$.
For each $i \notin D$, we pick $b_i$ to be a random element,
chosen uniformly from the set of elements not yet seen -- 
that is, from $\A^{\ell} \setminus \bigcup_i b_i$. 

Define random variables $R_j$ 
to be the time waited between time $k$ and the
first appearance of value $b_j$.
$$ R_j = \min \{ t \geq 1: X_{k+t} = b_j \}, \mbox{ for $j=1, \ldots, k$.}$$
\end{definition}
\begin{lemma} \label{lem:roomcorr} 
Suppose that $(Z_i)$ is an independent identically distributed
finite alphabet process with entropy $H$.
If, as block length $\ell \tends$, the number of blocks $k(\ell) \tends$ in 
such a way that 
$\lim_{\ell \tends} k(\ell)^{3/2} \ell q_{\max}^{\ell} = 0$, then
the difference term
$$ \frac{\sum_{i=1}^k ( \log R_i - \log S_i)}{\sqrt{k}}$$
tends to zero in probability.
\end{lemma}
\begin{proof} The key observation is that $S_i = R_i + (k-i)$, unless $i \in
D^c$. Further, $1 \leq S_i \leq R_i + (k-i)$.
This means that we can decompose
\begin{eqnarray*}
\pr \left( \frac{ \sum_{i=1}^k \log S_i - \log R_i}{\sqrt{k}} \geq \delta
\right)
& \leq &
\pr \left( \frac{ \sum_{i=1}^k \log (1 + (k-i)/R_i)}{\sqrt{k}} \geq \delta
\right) \\
& \leq & \pr \left( \sum_{i=1}^k k^{1/2} \frac{1}{R_i} \geq \delta \right) \\
& \leq & \frac{\sqrt{k}}{\delta} \sum_{i=1}^k \ep \frac{1}{R_i} 
= O( k(\ell)^{3/2} \ell q_{\max}^{\ell}),
\end{eqnarray*}
since $R \sim \geo{p}$ has $\ep 1/R = -p \log p/(1-p) =O( q_{\max}^\ell
\ell)$.

Similarly, we know that
\begin{eqnarray*}
\pr \left( \frac{ \sum_{i=1}^k \log R_i - \log S_i}{\sqrt{k}} \geq \delta
\right)
& \leq &
\pr \left( \frac{ \sum_{i=1}^k \log R_i \I(i \in D^c)}{\sqrt{k}} \geq \delta
\right) \\
& \leq & \frac{1}{\delta \sqrt{k}} \left( \sum_{i=1}^k \ep (\log R_i)
\pr(i \in D^c) \right) \\
& = & O( k(\ell)^{3/2} \ell q_{\max}^{\ell}),
\end{eqnarray*}
since $ \pr(i \in D^c) = 1 - \pr(i \in D) \leq
\left( 1 - \left(1- q_{\max}^{\ell} \right)^k \right) \leq k q_{\max}^{\ell},$
independently of $R_i$, and $\ep \log R_i = O(\ell)$.
\end{proof}
\subsection{Mean and variance of $\log R_i$}
We first find the leading order terms in $\ep \log R_i$ and $\var
\log R_i$ for all $i$. 
We use the following Lemma, the simplest form of the Euler-Maclaurin sum
formula.
\begin{lemma} \label{lem:eulmac}
For any differentiable function $f$ such that $f(x) \rightarrow 0$ 
as $x \tends$,
\begin{equation} \label{eq:eulmac}
\left| \sum_{i=1}^\infty f(i) - \int_1^\infty f(x) dx \right| 
\leq \frac{1}{2} |f(1)| + \frac{1}{2} \int_1^\infty |f'(x) | dx.
\end{equation}
\end{lemma}
\begin{proof}
Note that (integrating by parts) for all differentiable functions $f$, and
for all $a$:
$$ \int_{a}^{a+1} f(x) dx = \frac{1}{2} \left( f(a) + f(a+1) \right) 
- \int_a^{a+1} f'(x) \left( x- a- \frac{1}{2} \right),$$
which implies that
$$ \int_a^{a+1} f(x) dx - \frac{1}{2} \left( f(a) + f(a+1) \right) = R,$$
where $|R| \leq (\int_a^{a+1} |f'(x)| dx)/2$.
Summing such results from $a=1$ to $a=\infty$, we deduce that Equation 
(\ref{eq:eulmac}) holds. \end{proof}
\begin{lemma} \label{lem:values}
For $R \sim \geo{p}$:
\begin{eqnarray}
\mu(p) := \ep (\log R) 
& = & -\gamma - \log p + O(p) \label{eq:logrmean} \\
\sigma^2(p) := \var (\log R) & =  &
\frac{\pi^2}{6} + O(p \log p) \label{eq:logrvar} \\
 \ep | \log R - \mu(p) |^3 & \leq & K, \label{eq:logrmom}
\end{eqnarray}
where, as before, $\gamma$ is Euler's constant, and $K$ is a finite constant.
\end{lemma}
\begin{proof} Take $c = -\log(1-p)$ and
$f(x) = e^{-cx} \log x$, and use the fact that $|f'(x)| 
= | e^{-cx}/x - c e^{-cx} \log x| \leq | e^{-cx}/x| + 
c| e^{-cx} \log x|$.
We deduce by Equation (\ref{eq:eulmac}) that the difference
between the integral $I$ and sum $S$:
\begin{eqnarray*} 
\left| S - I \right| & = &
\left| \sum_{i=1}^\infty e^{-ci} \log i - \int_{1}^\infty
e^{-cx} \log x dx \right|  \\
& \leq & 
\frac{1}{2} \int_1^\infty \frac{e^{-cx}}{x} dx
+ \frac{c}{2} \int_1^\infty e^{-cx} \log x dx \\
& = & \left[
\frac{1}{2} e^{-cx} \log x \right]^{\infty}_1
+  c \int_1^\infty e^{-cx} \log x dx \\
& = & c I. 
\end{eqnarray*}
Using the fact (writing $\Gamma(0,\cdot)$ for the incomplete gamma
function) that $I = \Gamma(0,c)/c = (-\gamma - \log c + c + O(c^2))$,
we deduce that
$$ \ep \log R = \sum_{i=1}^\infty p(1-p)^{x-1} \log x 
= \frac{p}{1-p} \frac{\Gamma(0,c)}{c} (1 + \epsilon)
= (-\gamma - \log p) + O(p).$$
In the same way, for
$f(x) = e^{-cx} (\log x)^2$, the derivative $|f'(x)| 
= | 2e^{-cx} \log x/x - c e^{-cx}(\log x)^2| \leq | 2e^{-cx} \log x/x| + 
c| e^{-cx} (\log x)^2|$. As above
\begin{eqnarray*} 
\left| S - I \right| & = &
\left| \sum_{i=1}^\infty e^{-ci} \log i - \int_{1}^\infty
e^{-cx} \log x dx \right|
\leq c I. 
\end{eqnarray*}
Using the fact that $I = \pi^2/(6c) + (-\gamma + \log c)^2/c +
1 - O(c)$,
we deduce that
$ \ep \log R^2 - \mu(p)^2 = \pi^2/6 + O(p \log p).$

Finally, we bound the centred absolute third moment $\ep | \log R - \mu(p)|^3$.
We partition the real line into 3 intervals; firstly the set 
$A_1 = \{ x: |\log x - \mu(p)| \leq 1 \}$, secondly
$A_2 = \{ x: \log x - \mu(p) \geq 1 \}$, and thirdly
$A_3 = \{ x: \log x - \mu(p) \leq -1 \}$. We define integrals $K_i
= \ep |\log R - \mu(p)|^3 \I( x \in A_i)$ for $i=1,2,3$.
Clearly $K_1 \leq 1$. By Chernoff's bound, for $t \geq 1$:
\begin{eqnarray*} \pr( \log R - \mu(p) \geq t) 
& \leq & \frac{\ep R^s}{\exp( s(t + \mu(p)))} \\
& \leq & \exp(-2 t) \frac{(2-p)/p^2}{\exp(-2 \gamma + O(p))/p^2}
\leq 2 \exp(2 \gamma) \exp(-2 t),\end{eqnarray*}
taking $s = 2$. Hence  $K_2 = \int_1^\infty 3t^2 \pr( \log R - \mu(p) \geq t) 
dt \leq 4$. In a similar fashion,
$$ \pr( \log R - \mu(p) \leq -t) = \pr( R \leq
\exp(-\gamma - t))  \leq (1 - \exp( - e^{-\gamma - t})),$$
so $K_3 \leq 4$. Overall then, we can take $K = 9$.
\end{proof}
\subsection{Asymptotic independence}
Next we prove a lemma that shows that the $R_i$ are approximately
independent, giving explicit bounds on the difference between the joint
probability distribution and the product of the marginals.
\begin{lemma} \label{lem:probbounds}
Suppose that $(Z_i)$ is an independent identically distributed
finite alphabet process with entropy $H$. 
For $(R_i)$ as defined in Definition \ref{def:definer}, for
any $s,m$, and for any $\vc{a} = (a_1, \ldots a_{m-1})$, writing
$\vc{R} = (R_1, \ldots R_{m-1})$:
$$ \left( 1 - \frac{p_m}{1-S^*} \right)^{s-1} \leq
\pr(R_m \geq s | \vc{R} = \vc{a})
\leq (1-p_m)^{s-m},$$
where $S^* = p_1 + \ldots + p_{m-1}$.
\end{lemma}
\begin{proof} Given the values of $R_1, \ldots R_{m-1}$, we can write down
an explicit expression for the distribution of $R_m$.
\begin{equation} \label{eq:prodform}
 \pr(R_m \geq s | \vc{R} = \vc{a})
= \prod_{i=1}^{s-1} \pr(X_{k+i} \neq b_m | \vc{R} = \vc{a})
\end{equation}
We consider this product (\ref{eq:prodform}) 
term by term, for each value of $i$.
If for some $j \leq m-1$, the $a_j = i$, then $X_{k+i} = b_j$, so 
automatically $X_{k+i} \neq b_m$, and so the contribution from that $i$ to the
product (\ref{eq:prodform}) is 1.

Otherwise, if $a_j \neq i$ for all $j \leq m-1$, then
$$ \pr(X_{k+i} \neq b_m | \vc{R} = \vc{a})
= 1 - \pr(X_{k+i} = b_m | \vc{R} = \vc{a})
= 1 - \frac{p_m}{1 - S_i},
$$
where $S_i = \sum_{j=1}^{m-1} p_j \I( a_j > i)$. This is a decreasing 
function in $S_i$. 

It is clear that 
the product (\ref{eq:prodform}) is maximised when the first $(m-1)$ values
of $a_i$ occur in the first $m-1$ places, that is when
$\{ a_1, \ldots, a_{m-1} \} = \{ 1, \ldots, m-1 \}$.
In this case the value of (\ref{eq:prodform}) becomes $(1-p_m)^{s-m}$.

Similarly,
the product (\ref{eq:prodform}) is minimised when $S_i$ is maximised for
each $i$, that is when $a_j \geq s$ for each $j$. In that case, 
$S_i = \sum_{j=1}^{m-1} p_j = S^*$ for each $i$, and the product becomes
$( 1 - p_m/(1-S^*))^{s-1}$.
\end{proof}
\begin{proposition} \label{prop:charcl}
Suppose that $(Z_i)$ is an independent identically distributed
finite alphabet process with entropy $H$. 
For $(R_i)$ as defined in Definition \ref{def:definer}, 
$$ \left| \ep \exp \left( \frac{i \theta}{\sqrt{k}} \sum_{j=1}^m 
\log R_j \right) -
\ep \exp \left( \frac{i \theta}{\sqrt{k}}  \sum_{j=1}^{m-1}
 \log R_i \right)
\ep \exp \left( \frac{i \theta}{\sqrt{k}} 
\log R_m \right) \right|$$ is $O(\ell k(\ell)^{1/2} q_{\max}^{\ell})$.
\end{proposition}
\begin{proof}
Adapting Equation (22) of Newman \cite{newman80}, for any complex
continuously differentiable 
functions $f$ and $g$, and for random variables $U$ and $V$:
\begin{equation} \label{eq:covexp}
 \cov( f(U), g(V)) = \int_{-\infty}^{\infty}  \int_{-\infty}^{\infty} 
f'(u) g'(v) H_{U,V}(u,v) du dv,\end{equation}
where  
\begin{eqnarray*} 
H_{U,V}(u,v) & = & \pr( U \geq u, V \geq v) - \pr(U \geq u) \pr(V \geq v) \\
& = & - \pr( U \geq u, V < v) + \pr(U \geq u) \pr(V < v)
\end{eqnarray*}
We take $U = \log R_m$ and $V = \sum_{i=1}^{m-1} \log R_i$.
We will find non-negative
functions $h_{-}$ and $h_{+}$ such that
$- h_{-}(u,v) \leq H_{U,V}(u,v) \leq h_{+}(u,v)$ for all $u$ and $v$.
Since $U$, $V$ take non-negative values only, if $u < 0$ or
$v< 0$ then $H_{U,V}(u,v) = 0$. Hence, 
taking $f(t) = g(t) = \exp(i \theta t/\sqrt{k})$,
Equation (\ref{eq:covexp}) simplifies to
\begin{eqnarray} 
\lefteqn{\left| \cov \left( \exp \left( \frac{i \theta U}{\sqrt{k}} \right), 
\;
\exp \left( \frac{i \theta V}{\sqrt{k}} \right) \right) \right|} 
\nonumber \\ 
& = & \left|  - \frac{\theta^2}{k} \int_{0}^{\infty}  \int_{0}^{\infty} 
\exp \left( \frac{i \theta u}{\sqrt{k}} \right) \exp
\left( \frac{i \theta v}{\sqrt{k}} \right) H_{U,V}(u,v) du dv \right| 
\nonumber \\
& \leq &
\left| \frac{\theta^2}{k} \right|  \int_{0}^{\infty}  \int_{0}^{\infty} 
|H_{U,V}(u,v)| dv du \nonumber \\
& \leq & \left| \frac{\theta^2}{k} \right|  
\left( \int_{0}^{\infty}  \int_{0}^{\infty} h_{-}(u,v) dv du +
\int_{0}^{\infty}  \int_{0}^{\infty} h_{+}(u,v) dv du \right).
\label{eq:newnew} 
\end{eqnarray}
We know that $\int_0^\infty \pr(U \geq u) = \ep U$, 
$\int_0^\infty \pr(V \geq v) =
\ep V$, and $\int_0^\mu \pr(V < v) dv + \int_{\mu}^\infty \pr(V \geq v) dv
= \ep |V - \mu|$. Further, we can evaluate
$$ f(p) = \int_{0}^\infty (1-p)^{e^u-1} du = 
\frac{ \Gamma(0,-\log(1-p))}{1-p}.$$
Since $f(p)$ has an increasing, but negative gradient, with
$-f'(p) \leq 1/(-(1-p) \log (1-p)) \leq 1/(p(1-p))$, we know that for
$p \leq q$:
$f(p) - f(q) \leq (q-p)/(p(1-p))$.
This means that 
\begin{eqnarray} 
\int_{0}^{\infty} (1-p_m)^{e^u-1} - \left( 1 - \frac{p_m}{1-S^*} 
\right)^{e^u-1} du 
= O(k(\ell) q_{\max}^{\ell}). 
\label{eq:firstterm}
\end{eqnarray}
We rearrange Lemma \ref{lem:probbounds} and sum over values of $\vc{a}$ 
such that $\sum \log a_j \geq v$ or $\sum \log a_j < v$. 
\begin{enumerate}
\item{
For $v \geq \ep V$, we find that
\begin{eqnarray*}
h_{-}(u,v) & \leq & \left( (1-p_m)^{e^u-1} -
\left( 1 - \frac{p_m}{1-S} \right)^{e^u -1}  \right) \pr(V \geq v) 
\end{eqnarray*}
and similarly can take $h_+(u,v) \leq (1-p_m)^{1-m} \pr(U \geq u) \pr(V 
\geq v)$. 

Thus, by (\ref{eq:firstterm}), over this region 
$\int h_-(u,v) du dv \leq O(k(\ell) q_{\max}^{\ell}) \ep |V-\mu| = 
O(k(\ell)^{3/2} q_{\max}^{\ell})$
and $\int h_+(u,v) du dv \leq (1-p_m)^{1-m} (\ep U) \ep |V-\mu| = 
O(\ell k(\ell)^{3/2} q_{\max}^{\ell})$. }
\item{
For $v \leq \ep V$, we find that
\begin{eqnarray*}
h_{+}(u,v) & \leq & \left( 
(1-p_m)^{e^u-1} -
\left( 1 - \frac{p_m}{1-S} \right)^{e^u -1}  \right) \pr(V \leq v) 
\end{eqnarray*}
and similarly can take $h_-(u,v) \leq (1-p_m)^{1-m} \pr(U \geq u) \pr(V 
\geq v)$.

As before the integrals satisfy
$\int h_+(u,v) du dv = O(k(\ell)^{3/2} q_{\max}^{\ell})$
and $\int h_-(u,v) du dv 
= O(\ell k(\ell)^{3/2} q_{\max}^{\ell})$.}
\end{enumerate}
Substituting in Equation (\ref{eq:newnew}), the result follows.
\end{proof}
\subsection{Completing the proof of Theorem \ref{thm:main}}
The Lyapunov Central Limit Theorem (see for example Theorem 4.9 of
\cite{petrov})
implies that for independent $Y_1, \ldots, Y_k$, where $Y_i$ has
mean $\mu_i$, variance $\sigma^2_i$ and finite centred absolute
third moment $m_i = \ep |Y_i - \mu_i|^3$, if
\begin{equation} \label{eq:lyapcon}
\frac{\sum_{i=1}^k m_i}{\left( \sum_{i=1}^k \sigma^2 \right)^{3/2}}
\rightarrow 0,\end{equation}
then
\begin{equation} \label{eq:lyap}
\frac{\sum_{i=1}^k (Y_i - \mu_i)}
{\sqrt{ \var( \sum_{i=1}^k Y_i )}} \cond N(0,1)
\end{equation}

\begin{proof}{\bf of Theorem \ref{thm:main}}
Define a sequence of independent random variables $(T_i)$ with
$T_i \sim R_i$ then
Lemma \ref{lem:values} (in particular Equations (\ref{eq:logrvar})
and (\ref{eq:logrmom})) shows that the Lyapunov condition (\ref{eq:lyapcon})
holds for $\log T_i$, so that given values $p_1, \ldots, p_k$,
\begin{equation} \label{eq:tclt} \frac{ \sum_{i=1}^k 
( \log T_j - \mu(p_j))}{\sqrt{k}} \cond N(0,v),\end{equation}
where $v = \frac{1}{k} \left( \sum_{i=1}^k \var (\log R_i)^2 \right)
= \pi^2/6 + O(p \log p)$, and (by the Law of Large Numbers)
$\sum_{i=1}^k \mu(p_i) \rightarrow k(-\gamma + H \ell \log 2)$.

By repeated use of Proposition \ref{prop:charcl} then:
$$ \left| \ep \exp \left( \frac{i \theta}{\sqrt{k}} \sum_{j=1}^k 
(\log R_j - \mu(p_j)) \right) -
\prod_{j=1}^k \ep \exp \left( \frac{i \theta}{\sqrt{k}}  
(\log R_j - \mu(p_j)) \right)
\right|
$$ is $O(\ell k(\ell)^{3/2} q_{\max}^{\ell})$, so if this quantity tends
to zero, then the Central Limit Theorem for $\log T_i$, Equation 
(\ref{eq:tclt}) carries over to
give a Central Limit Theorem for $\log R_i$.
\end{proof}
\section{Equidistribution and negative association} \label{sec:negass}
In the case of equidistributed random variables,
we can establish a Central Limit Theorem 
under weaker conditions on $k(\ell)$, using 
negative association. This property captures the sense of 
dependence in which  one random variable being  large forces the others to
be smaller. Formally:
\begin{definition} \label{def:negass}
A collection of real-valued 
random variables $(U_k)$ is negatively associated (NA) if, 
for all increasing functions $f_1$ and $f_2$, the covariance 
\begin{equation} \label{eq:negass}
\cov( f_1(U_i: i \in A_1), f_2(U_j: j \in A_2)) \leq 0,\end{equation}
where $f_1$ and $f_2$ take arguments in disjoint sets
of indices $A_1$ and $A_2$. \end{definition}
The negative association property proves useful in many situations, not least
because Newman \cite{newman80} shows that the Central
Limit Theorem holds for NA sequences of random variables. 
Further, if $(U_k)$ forms an NA sequence, then for any increasing function
$f$, the $(f(U_k))$ are also an NA sequence. 
\begin{proposition} \label{prop:rna} 
Suppose that  $(Z_i)$ is  an independent equidistributed  process with
finite alphabet $\A$.
The $(R_i)$ introduced in Definition
\ref{def:definer} are negatively associated. \end{proposition}
\begin{proof}
Given the ordering $R_{\tau(1)} < R_{\tau(2)} < \ldots < R_{\tau(k)}$, the
actual values satisfy $R_{\tau(1)} \sim \geo{kp}$, $R_{\tau(2)} - R_{\tau(1)}
\sim \geo{(k-1)p}$ (independently) and so on. That is,
if we define $W_i$ independent with
$W_i \sim \geo{ (k+1-i)p}$, for $i=1, \ldots, k$, and define 
$U_j = \sum_{i=1}^j
W_i$, then $R_{\tau(i)} = U_i$, or $R_i = U_{\tau^{-1}(i)}$.

As in the proof of Theorem 3.4 of 
Hu \cite{hu2}, it suffices to show that 
Equation (\ref{eq:negass}) holds for symmetric functions
$f_1$ and $f_2$, with
$A_1 = \{ 1, \ldots , p \}$ and $A_2 = \{ p+1, \ldots , k \}$.
Given the functions $f_1$ and $f_2$, define
\begin{eqnarray*}
f_1^*(i_1, \ldots, i_p) & = & \ep \left( f_1( U_{i_1}, \ldots, U_{i_p}) 
\right), \\
f_2^*(i_{p+1}, \ldots, i_k) & = & \ep \left( f_2( U_{i_{p+1}}, \ldots, 
U_{i_k}) \right).
\end{eqnarray*}
If any index $i_l$ increases (for $l \in \{ 1, \ldots ,p \}$), then
(since $U_1 < U_2 < \ldots < U_k$) so does $U_{i_l}$, and hence (since
$f_1$ is increasing) so does $f_1^*$. That is, $f_1^*$ is an increasing
function of $\{ i_1, \ldots, i_p \}$. Similarly $f_2^*$ is also increasing.

For any permutation $\tau$, we define the increasing functions $$ g_1^*(\tau) =
f_1^*(\tau^{-1}(1), \ldots, \tau^{-1}(p)) \mbox{ and }
g_2^*(\tau) = f_2^*(\tau^{-1}(p+1), \ldots, \tau^{-1}(k)).$$
Theorem 2.11 of Joag-Dev and Proschan \cite{joag-dev} gives that
the uniform distribution on the set of permutations is negatively 
associated, so that 
\begin{equation} \label{eq:joagdev}
 \ep (g_1^*(\tau) g_2^*(\tau)) \leq \ep(g_1^*(\tau)) \ep(g_2^*(\tau)).
\end{equation}
Now $\ep g_1^*(\tau) = \ep f_1^*(\tau^{-1}(1), \ldots, \tau^{-1}(p)) =
\ep f_1(U_{\tau^{-1}(1)}, \ldots, U_{\tau^{-1}(p)}) \\ = \ep 
f_1(R_1, \ldots, R_p)$.
Similarly, $\ep g_2^*(\tau) = \ep f_2(R_{p+1}, \ldots R_{k})$ and
$\ep g_1^*(\tau) g_2^*(\tau) = \ep f_1(R_1, \ldots, R_p) f_2(R_{p+1}, \ldots
R_k)$, so Equation (\ref{eq:joagdev}) implies Equation (\ref{eq:negass})
as required.\end{proof}
\begin{lemma} \label{lem:roomcorr2} 
Suppose that  $(Z_i)$ is  an independent equidistributed  process with
finite alphabet $\A$.
If as block length $\ell \tends$, the number of
blocks  $k(\ell)  \tends$  in  such a  way  that  $k(\ell) \ell |\A|^{-\ell}
\rightarrow 0$, then the difference term
$$ \frac{\left( \sum_{i=1}^k 
( \log R_i - \log S_i) \right)}{\sqrt{k}}$$
tends to zero in probability.
\end{lemma}
\begin{proof} As before, $S_i = R_i + (k-i)$, unless $i \in
D^c$. Further, $1 \leq S_i \leq R_i + (k-i)$.
This means that we can decompose
\begin{eqnarray*}
 |\log R_i - \log S_i| & \leq &
|\log R_i - \log (R_i + (k-i))| + | \log(R_i + (k-i)) - \log S_i | \\
& \leq & \frac{k}{R_i} + \log(R_i + k) I(i \in D^c) 
\end{eqnarray*}
Since the $R_i$ are negatively associated, so are $-k/R_i$, so that by
Cauchy-Schwarz, writing $p=|\A|^{-\ell}$,  
\begin{eqnarray*}
\lefteqn{ \ep \left( \sum_{i=1}^k ( \log R_i - \log S_i ) \right)^2} \\
& \leq & 2 \sum_{i=1}^k k^2 \ep \frac{1}{R_i^2}
+ 2 \sum_{i \in D^c} \left( \ep \log (R_i)^2 + k^2 \ep \frac{1}{R_i^2}
\right) \\
& \leq & 2 k^3 (p^2 + O(p^3)) +  2 kp \left( 
O((-\log p)^2)+  k^2(p^2 + O(p^3))  \right) 
\end{eqnarray*}
This follows since $\ep 1/R^2_i = p {\rm Li}_2(p)/(1-p) = p^2 + O(p^3)$, 
where ${\rm Li}_2$ is the dilogarithm function, and since
for any $i$,
$$\pr(i \in D^c) \leq \pr(1 \in D^c) = 1 - \pr(1 \in D) = 1 - (1-p)^k 
\leq p k,$$
independently of $R_i$.
The lemma follows on dividing by $k$.
\end{proof}
\begin{lemma} \label{lem:covbound}
Suppose that  $(Z_i)$ is  an independent equidistributed  process with
finite alphabet $\A$.
For any $i \neq j$, the $R_i$ defined in Definition 
\ref{def:definer} satisfy $$|\cov(R_i, R_j)| = O \left( \ell  |\A|^{-\ell}
\right).$$
\end{lemma}
\begin{proof}
From the negative association proved in Proposition \ref{prop:rna}, we know
that the covariance is negative, so we need only bound it from below.
For any $x$, we know that, writing $p = |\A|^{-\ell}$.
$$ \pr(R_j = y | R_i = x) = \left\{ \begin{array}{ll}
\frac{p}{1-p} \left( \frac{ 1-2p}{1-p} \right)^{y-1} & \mbox{ for $y < x$,}\\
\left( \frac{1-2p}{1-p} \right)^{x-1} (1-p)^{y-x-1} p & \mbox{ for $y > x$,}
\end{array} \right.
$$
This means that
\begin{eqnarray}
\lefteqn{ \ep (\log R_j | R_i = x) } \nonumber \\
& = & 
\left( \sum_{y=1}^{x-1}
\frac{p}{1-p} \left( \frac{ 1-2p}{1-p} \right)^{y-1} \log y +
\sum_{y=x+1}^{\infty} p
\left( \frac{1-2p}{1-p} \right)^{x-1} (1-p)^{y-x-1} \log y \right) 
\nonumber \\
& = & 
\left( \sum_{y=1}^{\infty}
\frac{p}{1-p} \left( \frac{ 1-2p}{1-p} \right)^{y-1} \log y \right) 
- \frac{p}{1-p} \left( \frac{1-2p}{1-p} \right)^{x-1} 
\log x \nonumber \\
& & + \sum_{y=x+1}^{\infty} \left(
p \left( \frac{1-2p}{1-p} \right)^{x-1} (1-p)^{y-x-1}  
- \frac{p}{1-p} \left( \frac{ 1-2p}{1-p} \right)^{y-1} \right) \log y 
\label{eq:term3}
\end{eqnarray}
In Equation (\ref{eq:term3})
each summand is positive, so
we can replace $\log y \geq \log (y - x)$, and write $z = y-x$ to
obtain that Equation (\ref{eq:term3}) is greater than
\begin{eqnarray*} 
\lefteqn{ \sum_{z=1}^{\infty} \left( p
\left( \frac{1-2p}{1-p} \right)^{x-1} (1-p)^{z-1}  
- \frac{p}{1-p} \left( \frac{ 1-2p}{1-p} \right)^{z+x-1} \right) \log z } \\
& = & \left( \frac{1-2p}{1-p} \right)^{x-1} \left(
\sum_{z=1}^{\infty} \left( p (1-p)^{z-1}   - \frac{p}{1-p} \left(
\frac{ 1-2p}{1-p} \right)^z \right) \log z \right)  \\
& = & \left( \frac{1-2p}{1-p} \right)^{x-1} \left( \mu(p) - \frac{1-2p}{1-p}
\mu \left( \frac{p}{1-p} \right) \right).
\end{eqnarray*}
Overall then, using the notation of Lemma \ref{lem:values},
\begin{eqnarray*}
\lefteqn{\ep \log R_i \log R_j } \\
 & \geq & \sum_{x=1}^\infty p(1-p)^{x-1} \mu \left( \frac{p}{1-p} \right) 
- \frac{p}{2(1-p)} \sum_{x=1}^{\infty} 2p(1-2p)^{x-1} (\log x)^2  \\
& & + \frac{1}{2} \left( \sum_{x=1}^{\infty} 2p (1-2p)^{x-1} \log x \right)
\left( \mu(p) - \frac{1-2p}{1-p} \mu \left( \frac{p}{1-p} \right)  \right) \\
& = & \mu(p) \mu\left( \frac{p}{1-p} \right) 
- \frac{p}{2(1-p)} (\sigma^2(2p) + \mu(2p)^2) \\
& & + \frac{\mu(2p)}{2} \left( \mu(p) - \frac{1-2p}{1-p} \mu
\left( \frac{p}{1-p} \right)  \right).
\end{eqnarray*} 
Expanding this using Lemma \ref{lem:values}, we deduce that
$ \cov( \log R_i, \log R_j) \geq -p \log p$. Indeed, asymptotically,
$\cov( \log R_i, \log R_j) \geq p( (-\log p)/2 + \gamma)$. 
\end{proof}
We can now deduce the Central Limit Theorem for $\log S_i$:
\begin{proposition} \label{prop:sclt} 
Suppose that $(Z_i)$ is an independent equidistributed
finite alphabet process with entropy $H$.
If, as block length $\ell \tends$, the number of blocks $k(\ell) \tends$ in 
such a way that 
$k(\ell) \ell |\A|^{-\ell} \rightarrow 0$, then
$$ \frac{
\sum_{i=1}^{k(\ell)} (\log S_i - \ell H \log 2 + \gamma)}{\sqrt{
k(\ell) \pi^2/6 }} \cond N(0,1).$$
\end{proposition} 
\begin{proof} By Lemma \ref{lem:roomcorr2}, we need only prove the 
corresponding result for $\log R_i$.
By Proposition \ref{prop:rna}, the $R_i$ are 
negatively associated, and hence so are $\log R_i$.
Since $H(u,v) \leq 0$ for all $u$, $v$, adapting Equation (\ref{eq:newnew})
as in Newman \cite{newman80},
 gives the following result: if $U_1, \ldots, U_k$ are 
negatively associated then:
\begin{equation} \label{eq:newclt}
\left| \ep \exp \left( \frac{i \theta}{\sqrt{k}} \sum_{j=1}^k U_j \right) 
- \prod_{j=1}^k \ep \exp \left( \frac{i \theta U_j}{\sqrt{k}} \right)
\right| \leq \frac{\theta^2}{k} 
\sum_{i \neq j} | \cov(U_i,U_j)|. \end{equation}
This means that, taking $\mu = \ell H \log 2- \gamma$,
and $\varphi$ for the characteristic function of the $N(0,v)$:
\begin{eqnarray*}
\lefteqn{\left| \ep \exp \left( \frac{i \theta}{\sqrt{k}} 
\sum_{i=1}^k (\log R_i - 
\mu)  \right) - \varphi(\theta) \right| } \\
& \leq & 
\left| \ep \exp \left( \frac{ i \theta}{\sqrt{k}} 
\sum_{j=1}^k (\log R_j - \mu) \right) - 
\prod_{j=1}^k \ep \exp \left( \frac{ i \theta}{\sqrt{k}} 
\left(\log R_j - \mu \right) \right) \right| \\
& & +
\left| \prod_{j=1}^k \ep \exp 
\left( \frac{i \theta}{\sqrt{k}} \left( \log R_j - \mu \right) \right) - 
\varphi(\theta) 
\right|.
\end{eqnarray*}
Equation (\ref{eq:newclt}) bounds
the first term by $k \theta^2 |\cov(R_i,R_j)| = 
O(k(\ell) \ell \A^{-\ell})$, 
so we can control that term. We control the second term by using
the Lyapunov Central Limit Theorem, Equation (\ref{eq:lyap}).
\end{proof}
\section{Computational results} \label{sec:comp}
We present the results of some calculations, based on simulations with
random number generators, and on the decimal digits of well-known constants. 
In each case, we calculate the value of the statistic
$$ \frac{
\sum_{i=1}^{k(\ell)} (\log S_i - \ell H \log 2 + \gamma)}{\sqrt{
k(\ell) \pi^2/6  }}. $$
We present the results plotted as a quantile-quantile
plot using $R$ -- the line connects the upper and lower quartiles. 
If the distribution of the statistic
were exactly $N(0,1)$, we would see the majority of the points lying very
close to the line $y=x$.

\begin{figure}[ht!]
\begin{center}
\begin{psfrags}
    \psfrag{0}[t][t]{0}
    \psfrag{1}[t][t]{1}
    \psfrag{2}[t][t]{2}
   \psfrag{3}[t][t]{3}
    \psfrag{-1}[t][t]{-1}
    \psfrag{-2}[t][t]{-2}
   \psfrag{-3}[t][t]{-3}
\psfrag{Theoretical Quantiles}[c]{\footnotesize Theoretical Quantiles}
\psfrag{Sample Quantiles}[c]{\footnotesize Sample Quantiles}
\psfrag{Normal Q-Q Plot}{}
\epsfxsize=2.4in
\leavevmode\epsfbox{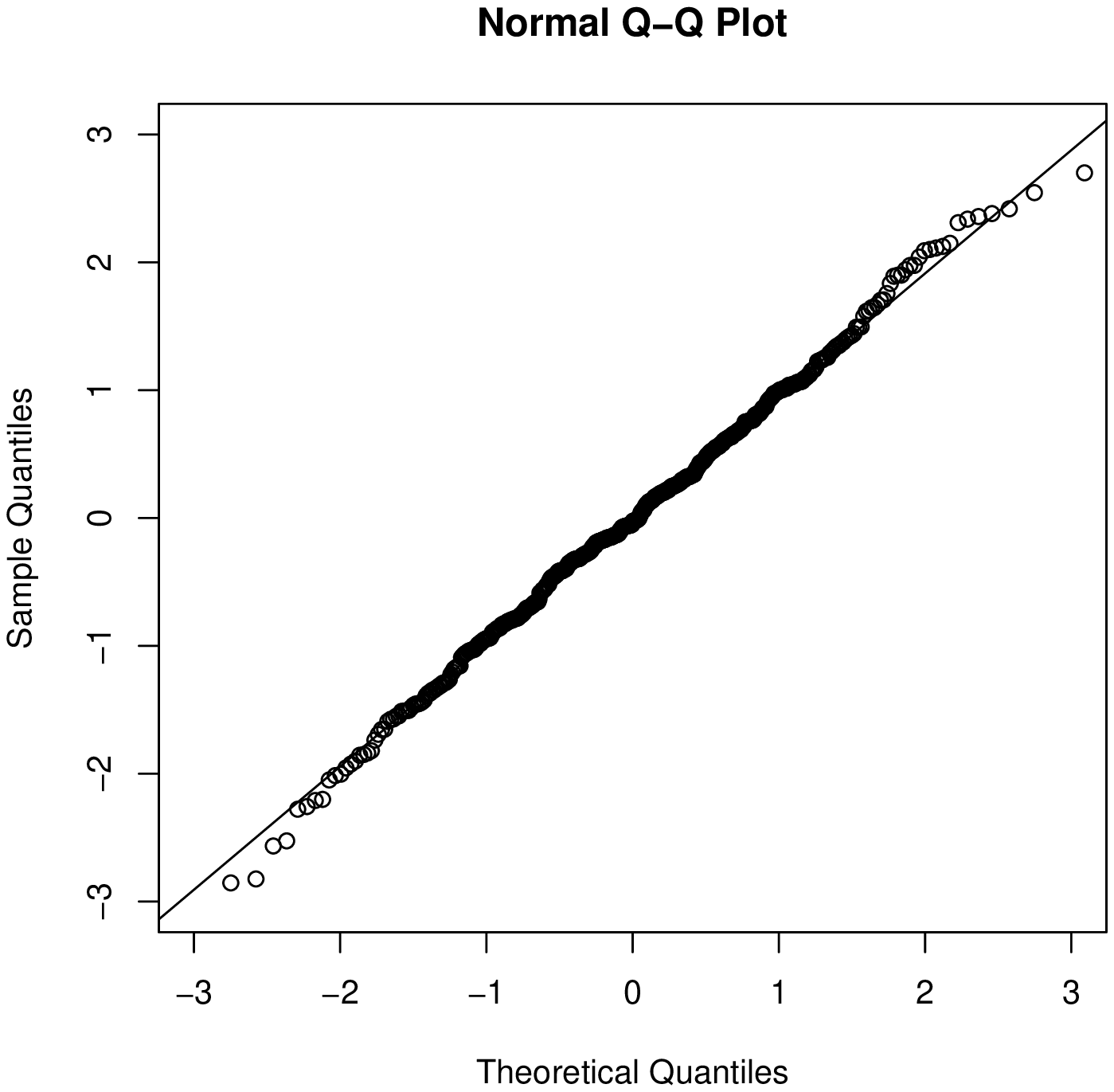}
\epsfxsize=2.4in
\leavevmode\epsfbox{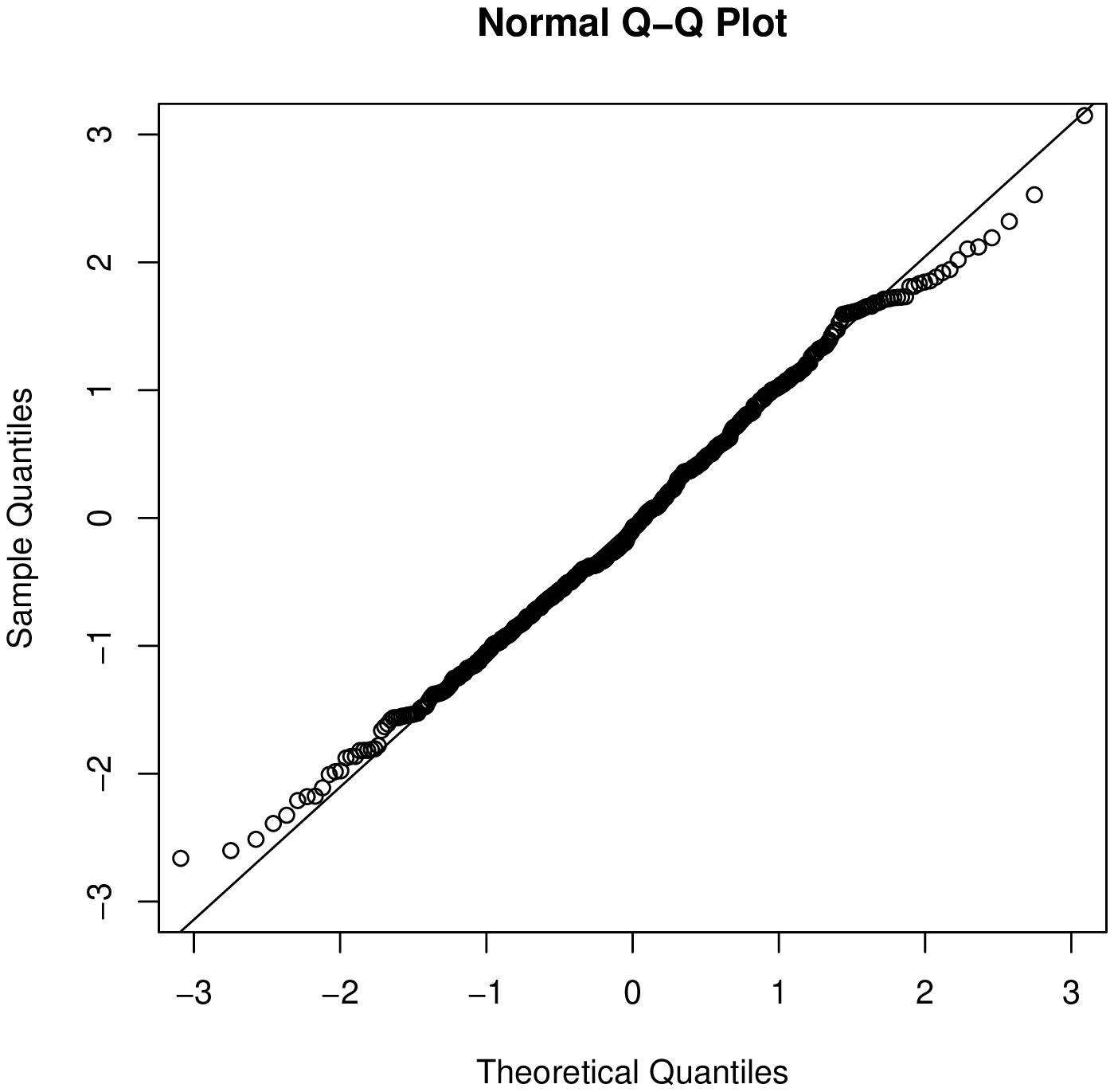}
\end{psfrags}
\caption{QQ plots: equidistributed binary data; 
(a) $k=250$, $\ell = 10$ (b) $k=1000$, $\ell =13$.}\label{fig:graph1} 
\end{center}\end{figure}

\begin{figure}[ht!]
\begin{center}
\begin{psfrags}
    \psfrag{0}[t][t]{0}
    \psfrag{1}[t][t]{1}
    \psfrag{2}[t][t]{2}
   \psfrag{3}[t][t]{3}
    \psfrag{-1}[t][t]{-1}
    \psfrag{-2}[t][t]{-2}
   \psfrag{-3}[t][t]{-3}
\psfrag{Theoretical Quantiles}[c]{\footnotesize Theoretical Quantiles}
\psfrag{Sample Quantiles}[c]{\footnotesize Sample Quantiles}
\psfrag{Normal Q-Q Plot}{}
\epsfxsize=2.4in
\leavevmode\epsfbox{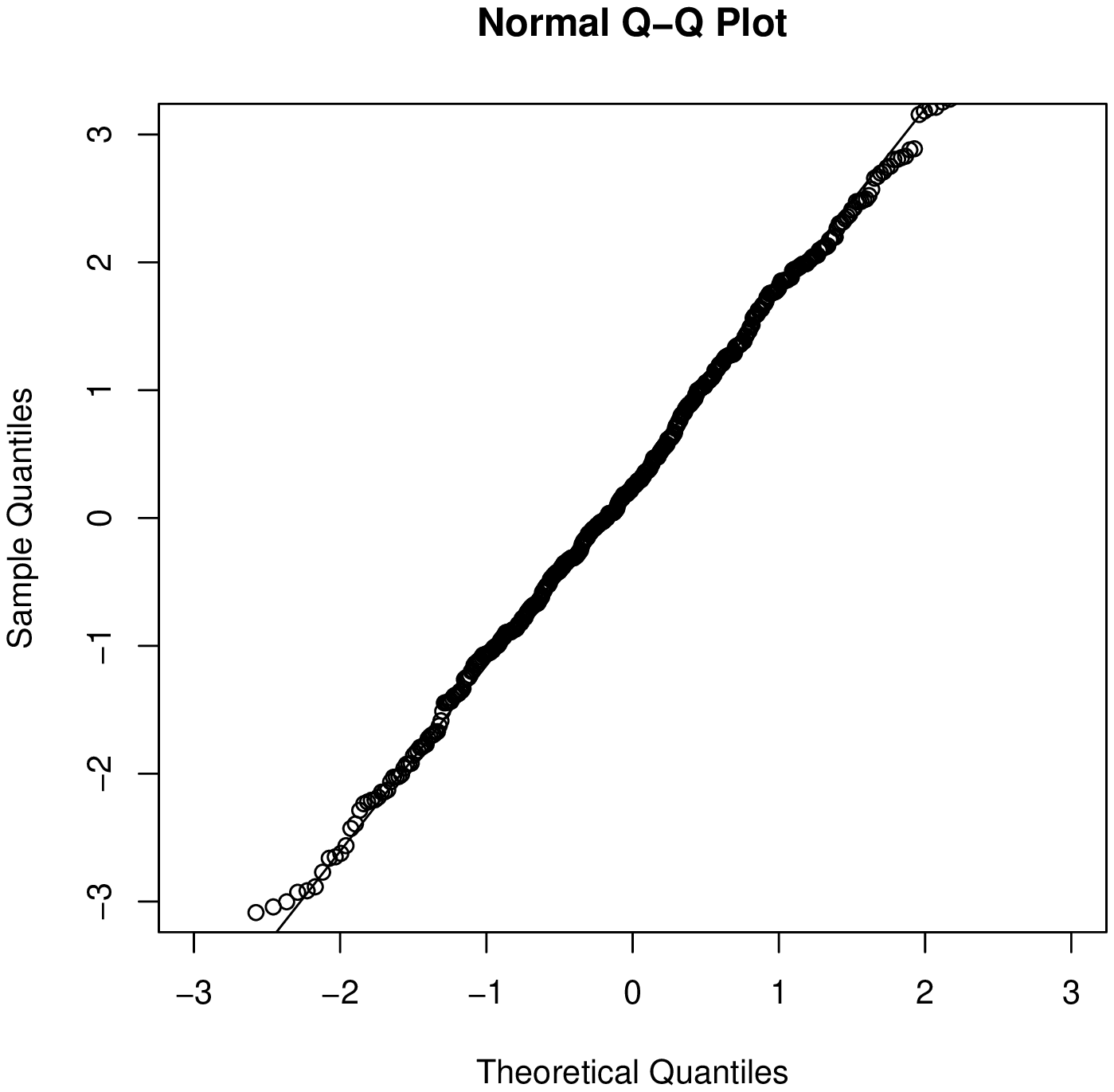}
\epsfxsize=2.4in
\leavevmode\epsfbox{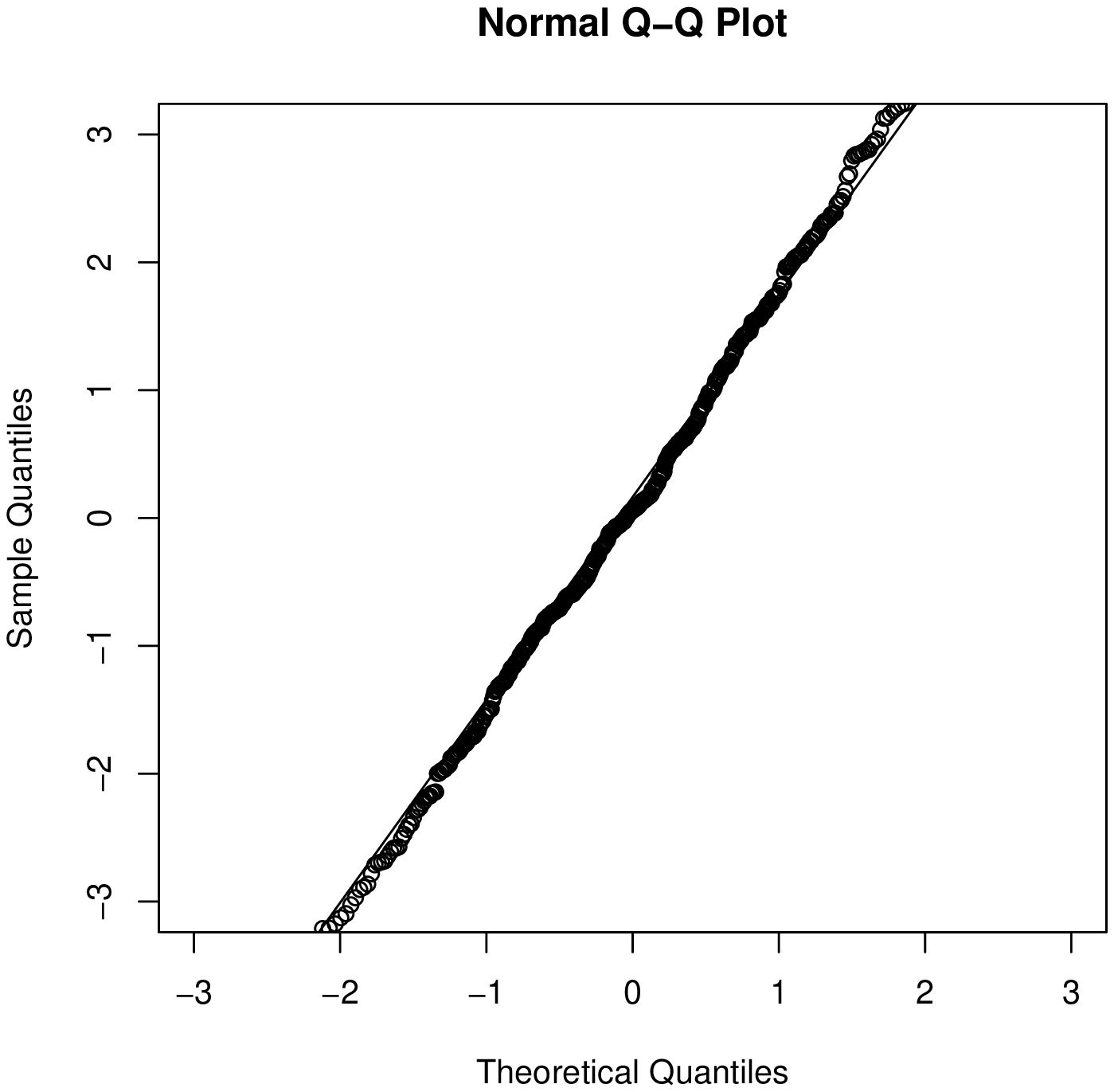}
\end{psfrags}
\caption{QQ plots: asymmetric ($\pr(Z_1 = 0) = 0.75$) binary data; 
(a) $k=250$, $\ell = 10$ (b) $k=1000$, $\ell =13$. }%
\label{fig:graph2} \end{center}\end{figure}

\begin{figure}[ht!]
\begin{center}
\begin{psfrags}
    \psfrag{0}[t][t]{0}
    \psfrag{1}[t][t]{1}
    \psfrag{2}[t][t]{2}
   \psfrag{3}[t][t]{3}
    \psfrag{-1}[t][t]{-1}
    \psfrag{-2}[t][t]{-2}
   \psfrag{-3}[t][t]{-3}
\psfrag{Theoretical Quantiles}[c]{\footnotesize Theoretical Quantiles}
\psfrag{Sample Quantiles}[c]{\footnotesize Sample Quantiles}
\psfrag{Normal Q-Q Plot}{}
\epsfxsize=2.4in
\leavevmode\epsfbox{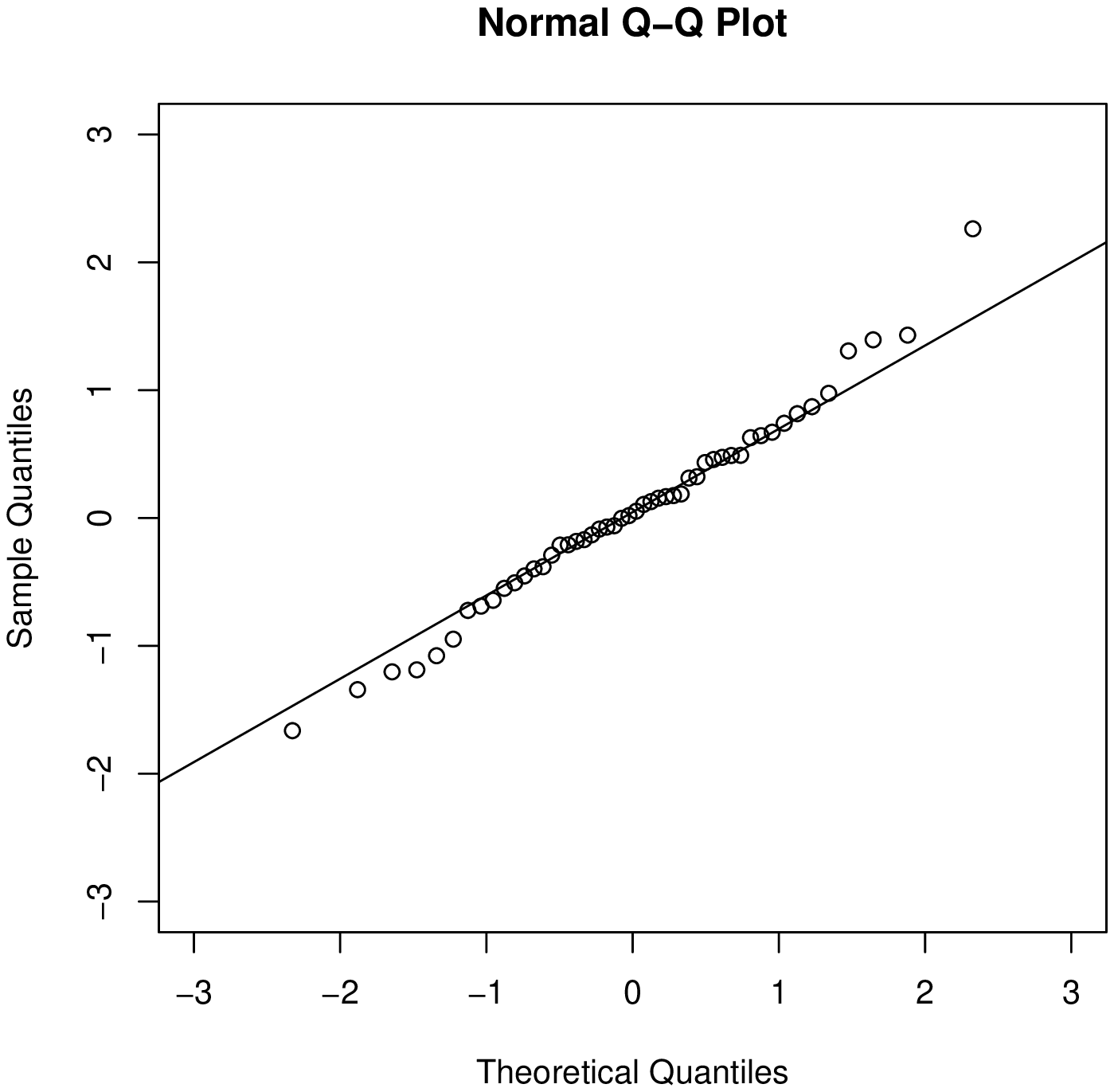}
\epsfxsize=2.4in
\leavevmode\epsfbox{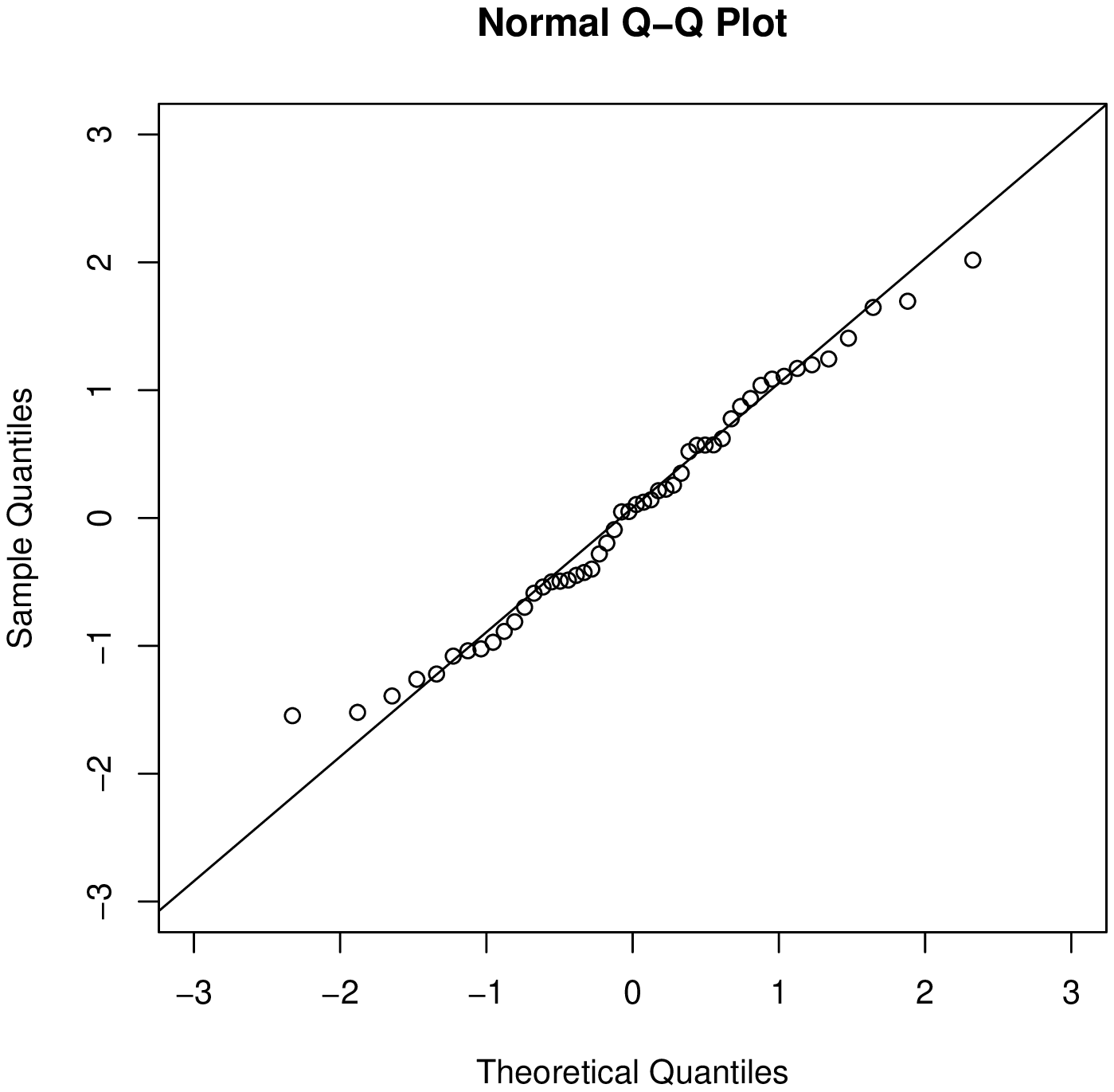}
\end{psfrags}
\caption{QQ plots: decimal data, $k=1000$, $\ell = 4$; (a) digits of $\pi$ 
(b) digits of $e$.}
\label{fig:graph3} \end{center}\end{figure}

To produce 
Figures \ref{fig:graph1} and \ref{fig:graph2}, we performed 500 trials on
simulated data. 
Figure \ref{fig:graph3} is based on breaking the first 
20 million decimal digits of $\pi$ and $e$ into 50 blocks of 400,000
digits each. We
used the program PiFast, which is freely downloadable from the Internet
\cite{gourdon}, and which can easily calculate tens of million digits of
constants such as $\pi$ and $e$.

In each case, the points do appear to lie on a straight line, though the 
sample variance is slightly smaller than expected. This could be remedied
by dividing by the square root of the true variance
$$ \var \left( \sum_{i=1}^{k(\ell)} \log S_i \right) =
k(\ell) \pi^2/6 + k(\ell) (k(\ell) -1) \cov(S_i, S_j) \leq 
k(\ell) \pi^2/6.$$
In order to do this, we would require an expansion, rather than simply an
approximation, for the covariance in Lemma \ref{lem:covbound} (since
the proof of Lemma \ref{lem:roomcorr2} shows that sums $\log R_i$ and 
$\log S_i$ have the same variance, asymptotically). Numerical
calculation suggests that $\cov(R_i, R_j) \sim (p \log p)/4$. Of course,
Lemma \ref{lem:covbound} only holds for equidistributed processes. However,
in general
the Asymptotic Equipartition Property suggests that we can assume that 
$\cov(R_i, R_j) \sim (2^{-H \ell} H \ell \log 2)/4$.

\clearpage
\bibliography{../bibliography/phd}
\end{document}